\documentclass[12pt,a4paper]{article}

\usepackage{latexsym,epsfig}
\usepackage{graphicx} 
\usepackage[usenames,dvipsnames]{color} 
\usepackage{url} 
\usepackage{soul} 
\usepackage{array} 
\usepackage{amssymb}
\usepackage{amsfonts}
\usepackage{amsmath}
\usepackage{hyperref}
\usepackage{verbatim}
\usepackage[sort&compress,square,comma,numbers]{natbib}
\usepackage{anysize}
\marginsize{1.5cm}{1.5cm}{1cm}{1cm}

\newtheorem{lemma}{Lemma}[section]
\newtheorem{proposition}{Proposition}[section]
\newtheorem{theorem}{Theorem}[section]
\newtheorem{corollary}{Corollary}[section]

\newtheorem{conjecture}{Conjecture}[section]
\newtheorem{observation}{Observation}[section]
\newcommand{\EndProof}{\hspace{\stretch{1}} $\Box$}

\newcommand{\pr}{\noindent{\bf Proof.}\ }

\newcommand{\G}{\Gamma}

\newcommand{\D}{\Delta}

\newcommand{\C}{\mathcal{C}}

\newcommand{\Jor}{J{\o}rgensen}
\newcommand{\rep}{\mathrm{rep}}

\newcommand{\M}{\mathrm{M}}
\newcommand{\N}{\mathrm{N}}
\newcommand{\g}{\mathrm{g}}
\newcommand{\m}{\mathrm{m}}

\title{On graphs of defect at most $2$}

\author{Ramiro Feria-Pur\'on$^{1,}$\footnote{\href{mailto:ramiro.feria.puron@gmail.com}{ramiro.feria.puron@gmail.com}}
\and
Mirka Miller$^{1,2,3,4,}$\footnote{\href{mailto:mirka.miller@newcastle.edu.au}{mirka.miller@newcastle.edu.au}}
\and
Guillermo Pineda-Villavicencio$^{5,}$\footnote{\href{mailto:work@guillermo.com.au}{work@guillermo.com.au} (Corresponding author)}\and
$\phantom{0}^{1}$\emph{\small School of
Electrical Engineering and Computer Science} \vspace{-1.5mm}\\
\emph{\small The University of Newcastle, Australia}\and
$\phantom{0}^{2}$\emph{\small Department of
Mathematics}\vspace{-1.5mm} \\\emph{\small University of West
Bohemia, Czech Republic}\and
$\phantom{0}^{3}$\emph{\small Department of Computer Science}\vspace{-1.5mm} \\\emph{\small King's College London, UK}\and
$\phantom{0}^{4}$\emph{\small Department of Mathematics}\vspace{-1.5mm} \\\emph{\small ITB Bandung, Indonesia} \and
$\phantom{0}^{5}$\emph{\small Centre for Informatics and Applied Optimization} \vspace{-1.5mm}\\
\emph{\small University of Ballarat, Australia}}

\begin{document}
\maketitle


\begin{abstract}
\noindent In this paper we consider the {\it degree/diameter problem}, namely, given natural numbers $\Delta\ge2$ and $D\ge1$, find the maximum number $\N(\Delta,D)$ of vertices  in a graph of maximum degree $\Delta$ and diameter $D$. In this context, the Moore bound $\M(\Delta,D)$ represents an upper bound for $\N(\Delta,D)$.

\noindent Graphs of maximum degree $\Delta$, diameter $D$ and order $\M(\Delta,D)$, called
{\it Moore graphs}, turned out to be very rare. Therefore, it is very interesting to investigate
graphs of maximum degree $\Delta\ge2$, diameter $D\ge1$ and order $\M(\Delta,D)-\epsilon$ with small
$\epsilon> 0$, that is, $(\Delta,D,-\epsilon)$-graphs. The parameter $\epsilon$ is called the
{\it defect}.

\noindent Graphs of defect $1$ exist only for $\Delta=2$.
When $\epsilon>1$, $(\Delta,D,-\epsilon)$-graphs represent a wide unexplored area.
This paper focuses on graphs of defect $2$. Building on the approaches developed in \cite{FP10} we obtain
several new important results on this family of graphs.

\noindent First, we prove that the girth of a $(\Delta,D,-2)$-graph with $\Delta\ge 4$ and $D\ge 4$ is $2D$.
 Second, and most important, we prove the non-existence of $(\Delta,D,-2)$-graphs with even
$\Delta\ge 4$ and $D\ge 4$; this outcome, together with a proof on the non-existence of
$(4,3,-2)$-graphs (also provided in the paper), allows us to complete the catalogue of
$(4,D,-\epsilon)$-graphs with $D\ge2$ and $0\le \epsilon\le2$. Such a catalogue is only the second census of $(\Delta,D,-2)$-graphs known at present, the first being the one of $(3,D,-\epsilon)$-graphs with $D\ge2$ and $0\le \epsilon\le2$ \cite{Jo92}.

\noindent Other results of this paper include necessary conditions for the existence of $(\Delta,D,-2)$-graphs with odd $\Delta\ge 5$ and $D\ge 4$, and the non-existence of $(\Delta,D,-2)$-graphs with odd $\Delta\ge 5$ and $D\ge 5$ such that $\Delta\equiv0,2\pmod{D}$.

\noindent Finally, we conjecture that there are no $(\Delta,D,-2)$-graphs with $\Delta\ge 4$ and $D\ge 4$, and comment on some implications of our results for the upper bounds of $\N(\Delta,D)$.
\end{abstract}

\noindent \textbf{Keywords:} Moore bound, Moore graph, Degree/diameter problem, defect, repeat.

\noindent\textbf{AMS  Subject Classification:} 05C35, 05C75.


\section{Introduction}
\label{sec:Introduction}

Due to the diverse features and applications of interconnection networks, it is possible to find many interpretations about network ``optimality'' in the literature. Here we are concerned with the following; see {\cite[pp. 168]{Hey96}}.

\begin{quote}
{\it An optimal network contains the maximum possible
number of nodes, given a limit on the number of connections
attached to a node and a limit on the distance between any two nodes of the network.}
\end{quote}

In graph-theoretical terms, this interpretation leads to the {\it degree/diameter
problem}, which can be stated as follows:

\begin{itemize}
\item[]{\it Degree/diameter problem}: Given natural numbers $\D\ge2$ and $D\ge1$, find the largest possible number $\N(\D,D)$ of vertices in a graph of maximum degree $\D$ and diameter $D$.
\end{itemize}

Note that $\N(\D,D)$ is well defined for $\D\ge2$ and $D\ge1$. An upper bound for
$\N(\D,D)$ is given by the {\it Moore bound} $\M(\D,D)$,
\[
\M(\D,D)=1+\D\left(1+(\D-1)+\dots+(\D-1)^{D-1}\right).
\]
Graphs of degree $\D$, diameter $D$ and  order $\M(\D,D)$ are called {\it Moore graphs}.

Only a few values of $\N(\D,D)$ are known at present. With the exception of
$\N(4,2)=\M(4,2)-2$ (see \cite{BA88}), $\N(5,2)=\M(5,2)-2$ (see \cite{NMb}),
$\N(6,2)=\M(6,2)-5$ (see \cite{Mol05}), $\N(3,3)=\M(3,3)-2$ (see \cite{Jo92}) and
$\N(3,4)=\M(3,4)-8$ (see \cite{Bus00}), the other known values of $\N(\D,D)$ are those for which
there exists a Moore graph.

Moore graphs are very rare. For $\D=2$ and $D\ge1$ they are the cycles on $2D+1$ vertices, whereas for
$D=1$ and $\D\ge2$ they are the complete graphs on $\D+1$ vertices. If $D=2$ and $\D\ge3$,
Moore graphs exist for $\D=3,7$ and possibly $57$, but not for any other degree; see \cite{HS60}. When $\D\ge3$ and $D\ge3$, there are no Moore graphs (\cite{Da73,BI73}).

Therefore, we are interested in studying the existence or otherwise of graphs of given maximum degree $\D$, diameter $D$ and order $\M(\D,D)-\epsilon$ for small $\epsilon>0$, that is, $(\D,D,-\epsilon)$-graphs, where the parameter $\epsilon$ is called the \emph{defect}.

The family of graphs of defect $\epsilon=1$ has been fully characterized; see \cite{EFH80,BI81,KT81}. For $\D=2$ and each $D\ge2$, the cycle on $2D$ vertices is the only $(2,D,-1)$-graph. For other values of $\D$ and $D$ there are no $(\D,D,-1)$-graphs.

Graphs of defect $\epsilon=2$ represent a wide unexplored area. The catalogue of $(3,D,-2)$ was
completed by \Jor~in \cite{Jo92}.  So far there have been several partial results achieved on
the existence or otherwise of $(\D,D,-2)$-graphs with $\D\ge4$ and $D\ge2$; see
\cite{El64,Jo92,BA88,NMb,MNP09,CG08} for $D=2$ and \cite{MS05,PM06} for $\D=4,5$. While the
paper \cite{MS05} claimed to have proved the non-existence of  $(4,D,-2)$-graphs for $D\ge3$, it
 turns out that the proof contained a mistake, so that only structural properties of $(4,D,-2)$-graphs were obtained. As a consequence, for  $(\D,D,-2)$-graphs with $\D\ge4$ and $D\ge2$ there has not been any definitive catalogue of any subfamily of such graphs until now.

For the sake of completeness we mention that, in the case of graphs with defect $\epsilon\ge3$, the only known work is the complete catalogue of $(3,D,-4)$-graphs provided in \cite{PM1}.

In this paper we consider $(\D,D,-2)$-graphs with $\D\ge 4$ and $D\ge 4$, and advance considerably
 the aforementioned question of the existence or otherwise of such graphs. To obtain our results we rely on combinatorial approaches which are inspired by those developed in \cite{FP10}.

Our first result is a proof that the girth of a $(\D,D,-2)$-graph with $\D\ge 4$ and $D\ge 4$
cannot be $2D - 1$ and therefore must be $2D$. Subsequently, we offer a non-existence proof of
$(\D,D,-2)$-graphs with even $\D\ge 4$ and $D\ge 4$. After ruling out the existence of
$(4,3,-2)$-graphs, we provide the first catalogue of $(\D,D,-\epsilon)$-graphs for $\D\ge4$,
$D\ge2$ and $0\le \epsilon\le2$, namely, the one of $(4,D,-\epsilon)$-graphs.

Other results of the paper include structural properties and necessary conditions for the
existence of $(\D,D,-2)$-graphs with odd $\D\ge 5$ and $D\ge 4$, and the non-existence of
$(\D,D,-2)$-graphs with odd $\D\ge 5$ and $D\ge 5$ such that $\D\equiv0,2\pmod{D}$.

Finally, we conjecture there are no $(\D,D,-2)$-graphs with $\D\ge 4$ and $D\ge 4$, and comment
 on the implications of our results for the upper bounds for $\N(\D,D)$.


\section{Known $(\D,D,-2)$-graphs}
\label{sec:KnownGraphs}

When $\D=2$ or $D=1$, there are no graphs of defect 2.

For $D=2$ there is a unique $(2,2,-2)$-graph (the path of length 2), exactly two non-isomorphic
$(3,2,-2)$-graphs, a unique $(4,2,-2)$-graph, and a unique $(5,2,-2)$-graph. All these graphs
are depicted in Fig.~\ref{fig:(Delta,D,-2)-graphs}. The non-existence
of $(\D,2,-2)$-graphs with $\D\ge6$, has been conjectured  but  not yet proved in spite of the partial support given in \cite{MNP09,CG08}.

\begin{figure}[!ht]
\begin{center}
\makebox[\textwidth][c]{\includegraphics[scale=1]{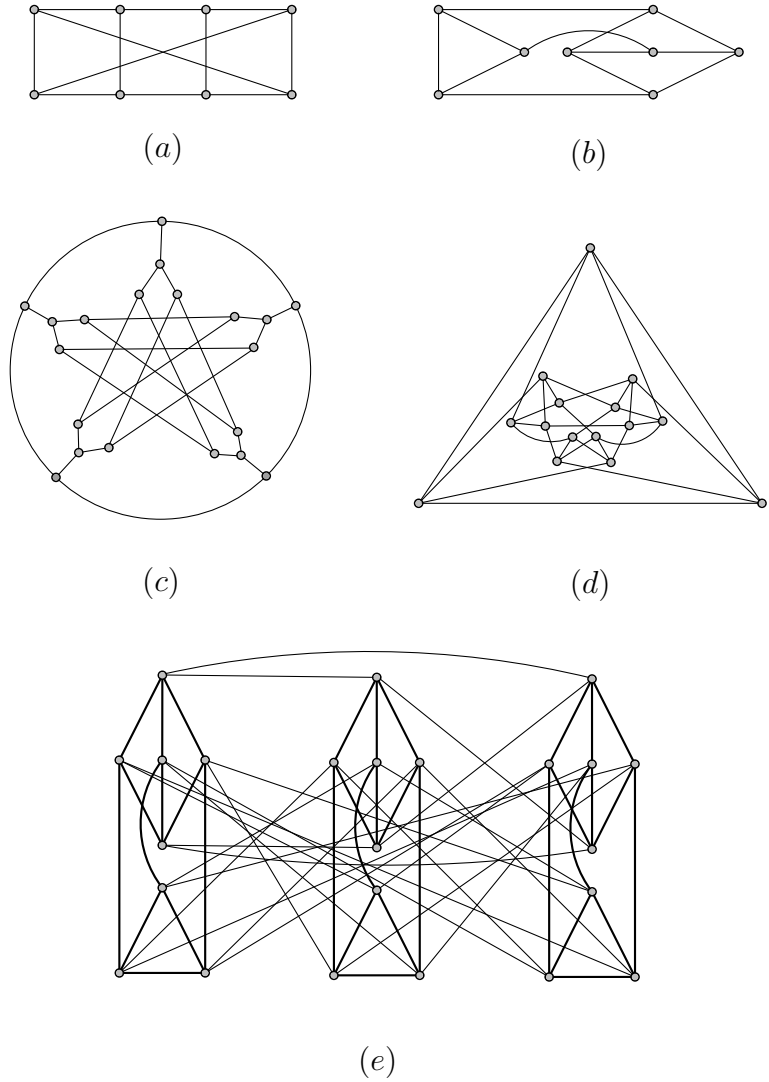}}
\caption{Known graphs of defect 2. $(a)$ and $(b)$ the two $(3, 2,-2)$-graphs, $(c)$ the unique
$(3, 3,-2)$-graph, $(d)$ the unique $(4, 2,-2)$-graph and $(e)$ the unique $(5, 2,-2)$ graph (note
that this graph is formed by connecting appropriately 3 copies of the graph $(b)$).}
\label{fig:(Delta,D,-2)-graphs}
\end{center}
\end{figure}

For $\D=3$ and $D\ge3$ there is a unique $(3,3,-2)$-graph, which is depicted in
Fig.~\ref{fig:(Delta,D,-2)-graphs} ($c$). This graph, together with the two aforementioned
$(3,2,-2)$-graphs, comprise the complete catalogue of cubic graphs of defect 2; see \cite{Jo92}.


\section{Notation and Terminology}
\label{sec:Notation}

The terminology and notation used in this paper is standard and consistent with that used in \cite{Die05}, so only those concepts that can vary from texts to texts will be defined.

All graphs considered in this paper are simple. The vertex set of a graph $\G$ is denoted by $V(\G)$, and its edge set by $E(\G)$. For an edge $e=\{x,y\}$, we write $x\sim y$. The set of neighbors of a vertex $x$ in $\G$ is denoted by $N(x)$.

A path of length $k$ is called a {\it $k$-path}. A path from a vertex $x$ to a vertex $y$ is
denoted by $x-y$. We use the following notation for subpaths of a path $P=x_0x_1\ldots x_k$:
$x_iPx_j = x_i\ldots x_j$, where $0\leq i\leq j\leq k$. A cycle of length $k$ is called a
{\it $k$-cycle}. The {\it girth} of $\G$, denoted g=g($\G$), is the length of the shortest cycle in $\G$.

The union of three independent paths of length $D$ with common endvertices is denoted by
$\Theta_D$. In a graph $\G$, a vertex of degree at least 3 is called a \emph{branch vertex} of
$\G$.


\section{Preliminary Results}
\label{sec:Preliminary}

We begin this section with a known condition for the regularity of a $(\Delta,D,-\epsilon)$-graph, which can be easily deduced by considering the existence of a vertex of degree at most $\D-1$ in
such a graph.

\begin{proposition}
\label{prop:regularity} For $\epsilon<1+(\Delta-1)+(\Delta-1)^2+\ldots+(\Delta-1)^{D-1}$, $\D\ge 3$ and $D\ge 2$, a $(\D,D,-\epsilon)$-graph is regular.
\end{proposition}

By Proposition \ref{prop:regularity}, a $(\D,D,-2)$-graph $\G$ with $\D\ge3$ and $D\ge2$ must be
 regular; we therefore use the symbol $d$ rather than $\D$ to denote the degree of $\G$, as is customary. We call a cycle of length at most $2D$ in $\G$ a {\it short cycle}.

\begin{proposition}[Lemma 2 from \cite{Jo92}]
\label{prop:Girth}
Let $\G$ be a $(d,D,-2)$-graph  with $d\ge3$ and $D\ge 2$. Then $2D-1\le \g(\G)\le 2D$. Furthermore, if $x$ is a vertex in $\G$ then either

\begin{itemize}
\item[$(i)$] $x$ is contained in one $(2D-1)$-cycle and no other short cycle; or

\item[$(ii)$] $x$ is contained in one $\Theta_D$, and every short cycle containing $x$ is contained in this $\Theta_D$; or

\item[$(iii)$] $x$ is contained in exactly two $2D$-cycles whose intersection is a $\ell$-path with $0\le\ell\le D-1$, and no other short cycle.
\end{itemize}
\end{proposition}

Each case is considered as a type. For instance, a vertex satisfying case $(i)$ is called a vertex of Type $(i)$.

While the statements of Proposition \ref{prop:Girth} and \cite[Lemma 2]{Jo92} slightly differ, both assertions are clearly equivalent. However, the statement of Proposition \ref{prop:Girth} is more consistent with the presentation of our results and allows us to make the following observation, which will be used implicitly throughout the paper.

\begin{observation}
\label{obs:Girth}
Let $\G$ be a $(d,D,-2)$-graph  with $d\ge3$ and $D\ge 2$, and $C$ a short cycle in $\G$. Then all vertices in $C$ are of the same type.
\end{observation}

In view of Proposition \ref{prop:Girth}, we define the following concepts:

We say that the vertex $x'$  is a {\it repeat} of the vertex $x$ with {\it multiplicity} $\m_x(x')$ $(1 \leq \m_x(x')\le2)$ if there are exactly $\m_x(x')+1$ different paths of length at most $D$ from $x$ to $x'$. For vertices $x$ and $x'$ lying on a short cycle $C$, we denote the vertex $x'$ by $\rep^C(x)$ if $x$ and $x'$ are repeats.

A vertex $x$ is called {\it saturated} if $x$ cannot belong to any further short cycle. If two $2D$-cycles $C^1$ and $C^2$ are non-disjoint, we say that $C^1$ and $C^2$ are \emph{neighbor cycles}.

From now on, whenever we refer to paths we mean shortest paths. As in \cite{FP10}, we extend the concept of repeat to paths. For a path $P=x-y$ of length at most $D-1$ contained in a $2D$-cycle $C$, we denote by $\rep^C(P)$ the path $P'\subset C$ defined as $\rep^C(x)-\rep^C(y)$. We say that $P'$ is the \emph{repeat} of $P$ in $C$ and vice versa, or simply that $P$ and $P'$ are \emph{repeats} in $C$.


\begin{lemma}[Odd Saturating Lemma]
\label{lem:OddSaturating}
Let $\G$ be a $(d,D,-2)$-graph with $d\ge 4$ and $D\ge 2$, and $\C$ a $(2D-1)$-cycle in $\G$. Let $\alpha$ be a vertex in $\C$ with repeat vertices $\alpha'_1,\alpha'_2$ in $\C$, $\gamma$ a neighbor of $\alpha$ not contained in $\C$, and $\mu_1,\mu_2,\ldots,\mu_{d-2}$ the neighbors of $\alpha'_2$ not contained in $\C$.

Then there is in $\G$ a vertex $\mu\in\{\mu_1,\mu_2,\ldots,\mu_{d-2}\}$ and a $2D$-cycle $\C^1$ such that $\gamma$ and $\mu$ are repeats in $\C^1$, and $\C\cap\C^1=\emptyset$.
\end{lemma}

\begin{figure}[!ht]
\begin{center}
\makebox[\textwidth][c]{\includegraphics[scale=1]{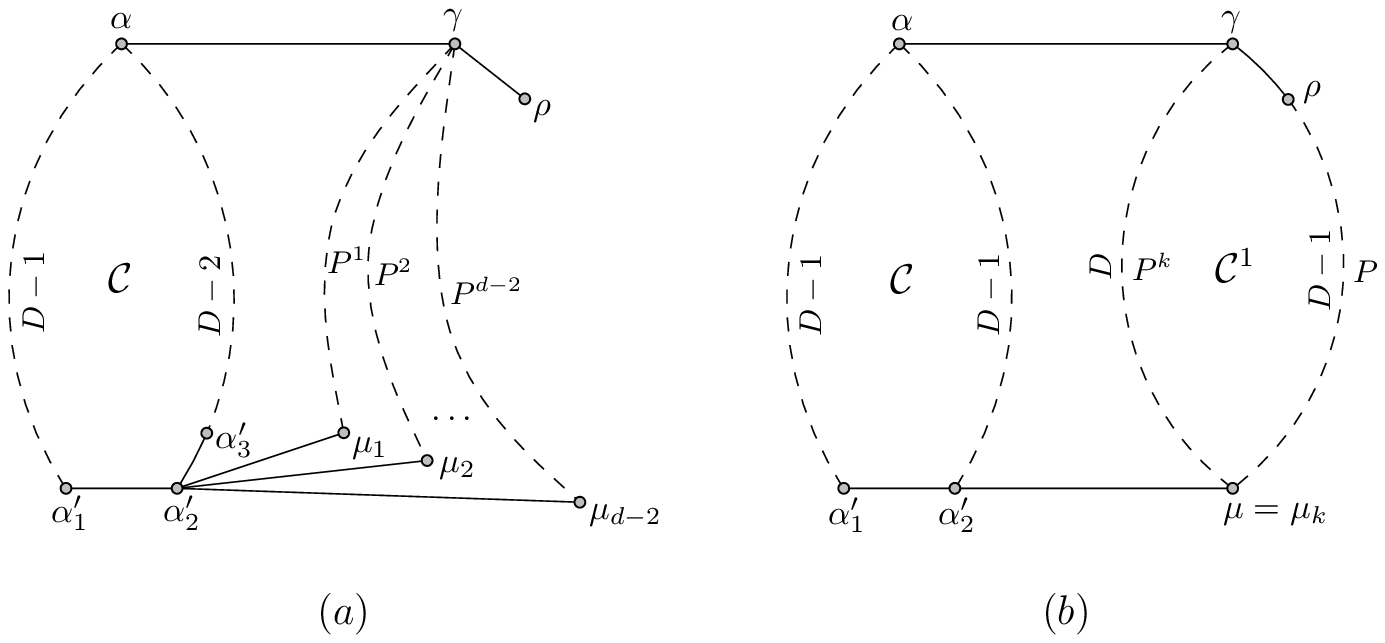}}
\caption{Auxiliary figure for Lemma \ref{lem:OddSaturating}}
\label{fig:OddSaturating}
\end{center}
\end{figure}

\pr Let $\alpha'_3$ be the neighbor of $\alpha'_2$ in $\C$ other than $\alpha'_1$. For $1\le i\le d-2$, consider the path $P^i=\gamma-\mu_i$. Since all vertices in $\C$ are saturated, $P^i$ cannot go through $\C$ and must be a $D$-path, so $P^i\cap\C=\emptyset$. Also, it follows that $V(P^i\cap P^j)=\{\gamma\}$ for any $1\le i<j\le d-2$; otherwise either $\g(\G)<2D-1$ or the vertex $\alpha'_2$ would belong to an additional short cycle, both contradictions to Proposition \ref{prop:Girth}. See Fig.~\ref{fig:OddSaturating} $(a)$.

Let $\rho$ be a neighbor of $\gamma$ other than $\alpha$, not contained in any of the paths $P^1,P^2,\ldots,P^{d-2}$ (there is exactly one such vertex). Consider a path $P=\rho-\alpha'_2$. $P$ cannot go through $\alpha'_3$; otherwise there would be a second short cycle $\rho P\alpha'_3\C\alpha\gamma\rho$ in $\G$ containing $\alpha$. Similarly, $P$ cannot go through $\alpha'_1$ and consequently, it must go through a vertex $\mu_k\in\{\mu_1,\mu_2,\ldots,\mu_{d-2}\}$. Finally note that, since all vertices in $\C$ are saturated and $2D-1\le \g(\G)\le 2D$, $P$ must be a $D$-path, $V(P\cap P^k)=\{\mu_k\}$ and $V(P\cap\C)=\{\alpha'_2\}$ .

This way, we obtain there is a vertex $\mu=\mu_k$ and a $2D$-cycle $\C^1=\gamma\rho P\mu P^k\gamma$ such that $\gamma$ and $\mu$ are repeats in $\C^1$, and $\C\cap\C^1=\emptyset$ (Fig.~\ref{fig:OddSaturating} $(b)$).\EndProof

\begin{lemma}[Even Saturating Lemma]
\label{lem:EvenSaturating}
Let $\G$ be a $(d,D,-2)$-graph with $d\ge 4$ and $D\ge 2$, and $\C$ a $2D$-cycle in $\G$. Let $\alpha,\alpha'$ be two vertices in $\C$ such that $\alpha'=\rep^{\C}(\alpha)$, $\gamma$ a neighbor of $\alpha$ not contained in $\C$, and $\mu_1,\mu_2,\ldots,\mu_{d-2}$ the neighbors of $\alpha'$ not contained in $\C$. Suppose there is no short cycle in $\G$ containing the edge $\alpha\sim\gamma$ and intersecting $\C$ at a path of length greater than $D-2$.

Then there is in $\G$ a vertex $\mu\in\{\mu_1,\mu_2,\ldots,\mu_{d-2}\}$ and a short cycle $\C^1$ such that $\gamma$ and $\mu$ are repeats in $\C^1$, and $\C\cap\C^1=\emptyset$.
\end{lemma}

\begin{figure}[!ht]
\begin{center}
\makebox[\textwidth][c]{\includegraphics[scale=1]{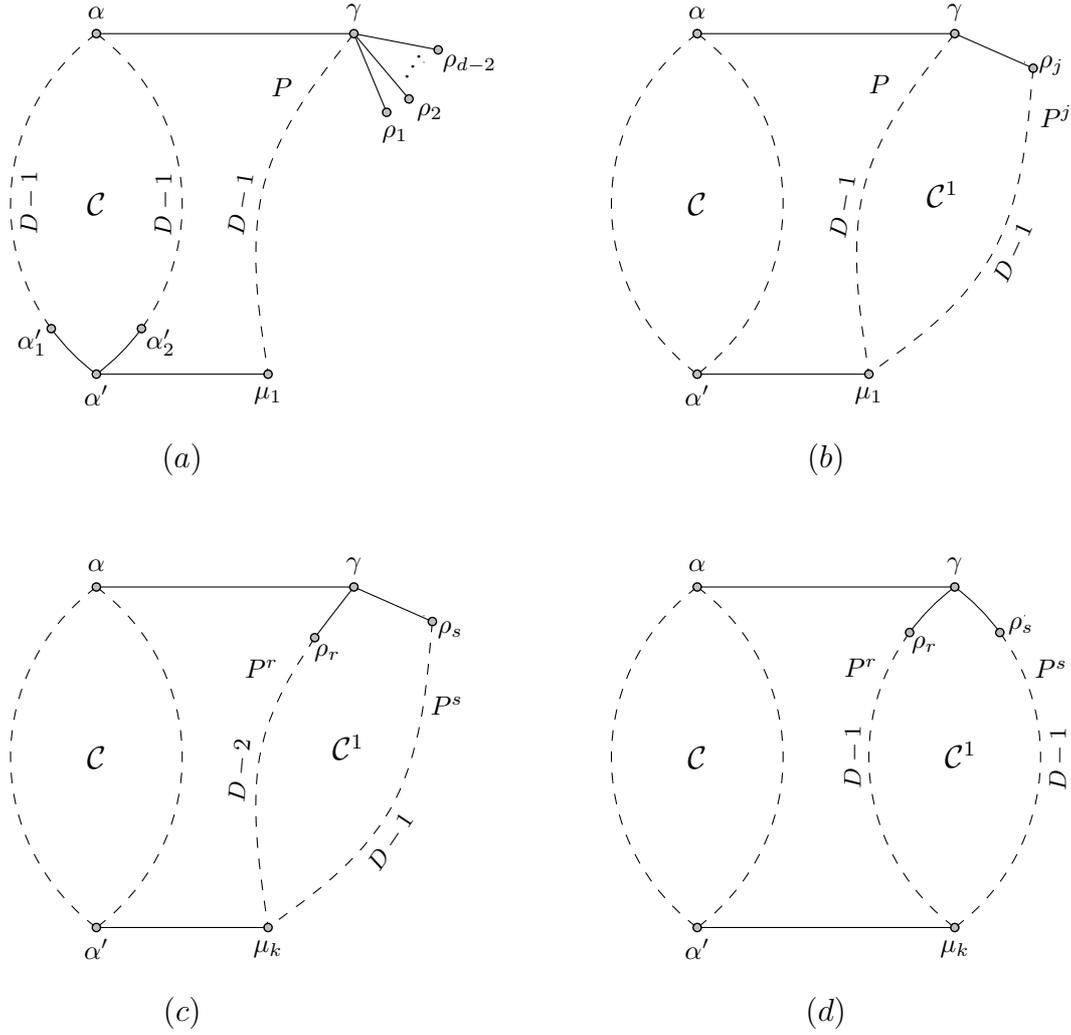}}
\caption{Auxiliary figure for Lemma \ref{lem:EvenSaturating}}
\label{fig:EvenSaturating}
\end{center}
\end{figure}

\pr Let $\alpha'_1,\alpha'_2$ be the neighbors of $\alpha'$ contained in $\C$. First, consider a path $P=\gamma-\alpha'$. Since there is no short cycle in $\G$ containing the edge $\alpha\sim\gamma$ and intersecting $\C$ at a path of length greater than $D-2$, $P$ must be a $D$-path and cannot go through $\alpha'_1$ or $\alpha'_2$. Therefore, the path $P$ must go through one of the neighbors of $\alpha$ not contained in $\C$ (say $\mu_1$). In addition, we have that $V(P\cap\C)=\{\alpha'\}$. See Fig.~\ref{fig:EvenSaturating} $(a)$.

Let $\rho_1,\rho_2,\ldots,\rho_{d-2}$ be the neighbors of $\gamma$ other than $\alpha$, not
contained in $P$. For $1\le i\le d-2$, consider the path $P^i=\rho_i-\alpha'$. Since there is
no short cycle in $\G$ containing the edge $\alpha\sim\gamma$ and intersecting $\C$ at a path of
 length greater than $D-2$, $P^i$ must have length at least $D-1$ and cannot contain any of the
vertices in $\{\alpha'_1,\alpha'_2,\gamma\}$. Consequently, $P^i$ must go through one of the vertices in $\{\mu_1,\mu_2,\ldots,\mu_{d-2}\}$. Note also that $V(P^i\cap\C)=\{\alpha'\}$ and that $V(P^i\cap P^j)\subseteq\{\alpha'\}\cup\{\mu_1,\mu_2,\ldots,\mu_{d-2}\}$, for any $1\le i<j\le d-2$.

If, for some $j$ $(1\le j\le d-2)$, the path $P^j$ goes through $\mu_1$ then $P^j$ must be a $D$-path and there is a $(2D-1)$-cycle $\C^1=\gamma P\mu_1P^j\rho_j\gamma$ in $\G$ such that $\gamma$ and $\mu=\mu_1$ are repeats in $\C^1$, and $\C\cap\C^1=\emptyset$. This case is depicted in Fig.~\ref{fig:EvenSaturating} $(b)$.

If, on the other hand, there is no $j$ $(1\le j\le d-2)$ such that $P^j$ goes through $\mu_1$
 then there must exist a vertex $\mu_k$ $(2\le k\le d-2)$ and paths $P^r,P^s$
$(1\le r<s\le d-2)$ such that both $P^r$ and $P^s$ go through $\mu_k$. Since $\g(\G)\ge 2D-1$,
at most one of the paths $P^r,P^s$ has length $D-1$. If one of these paths (say $P^r$) has
length $D-1$ then there is a $(2D-1)$-cycle $\C^1=\gamma\rho_r P^r\mu_kP^s\rho_s\gamma$ in
$\G$ such that $\gamma$ and $\mu=\mu_k$ are repeats in $\C^1$, and $\C\cap\C^1=\emptyset$ (as in
 Fig.~\ref{fig:EvenSaturating} $(c)$). If both $P^r$ and $P^s$ are $D$-paths then there is a $2D$-cycle $\C^1=\gamma\rho_r P^r\mu_kP^s\rho_s\gamma$ in $\G$ such that $\gamma$ and $\mu=\mu_k$ are repeats in $\C^1$, and $\C\cap\C^1=\emptyset$ (as in Fig.~\ref{fig:EvenSaturating} $(d)$).\EndProof

\subsection{Repeats of Cycles}

The extension of the concept of repeat to short cycles was introduced in \cite{FP10} in the context of bipartite graphs missing the bipartite Moore bound by 4 vertices. Here, inspired by the ideas put forward in \cite{FP10}, we extend the concept of repeat to $2D$-cycles of graphs of defect 2; see the Repeat Cycle Lemma.

\begin{lemma}[Repeat Cycle Lemma]
\label{lem:RepeatCycle}
Let $\G$ be a $(d,D,-2)$-graph with $d\ge 4$ and $D\ge 2$, and $C$ a $2D$-cycle in $\G$. Let
$\{C^1,C^2,\ldots, C^k\}$ be the set of neighbor cycles of $C$, and  $I_i=C^i\cap C$ for
$1\le i\le k$. Suppose at least one $I_j$, for $j\in\{1,\ldots, k\}$, is a path of length
smaller than $D-1$. Then there is an additional $2D$-cycle $C'$ in $\G$ intersecting $C^i$ at $I'_i=\rep^{C^i}(I_i)$, where $1\le i\le k$.
\end{lemma}

\begin{figure}[!ht]
\begin{center}
\makebox[\textwidth][c]{\includegraphics[scale=.9]{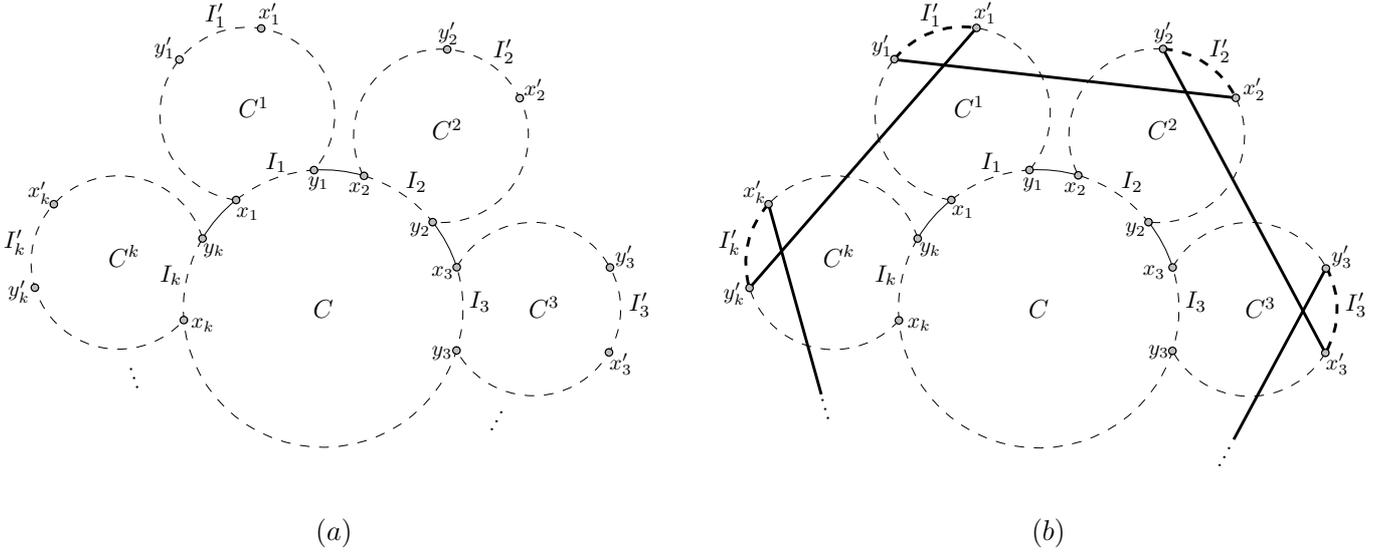}}
\caption{Auxiliary figure for Lemma \ref{lem:RepeatCycle}}
\label{fig:RepeatCycle}
\end{center}
\end{figure}

\pr We denote the neighbors of $C$ by $C^1,C^2,\ldots C^k$ and their corresponding intersection
paths with $C$ by $I_1=x_1-y_1,I_2=x_2-y_2,\ldots,I_k=x_k-y_k$ in such a way that $C=x_1I_1y_1x_2I_2y_2\ldots x_kI_ky_kx_1$. For $1\le i\le k$, we also denote the repeats of $I_i$ by $I'_i=x'_i-y'_i$, where $x'_i=\rep^{C^i}(x_i)$ and $y'_i=\rep^{C^i}(y_i)$ (see Fig.~\ref{fig:RepeatCycle} ($a$)).

For $1\le i \le k$, consider the cycles $C^i$ and $C^{(i \bmod{k})+1}$.

Suppose that $I_i$ is a path of length smaller than $D-1$. Since $y_i$ is saturated, there cannot be a short cycle in $\G$, other than $C$, containing the edge $y_i\sim x_{(i \bmod{k})+1}$. Since $I_i$ is a path of length smaller than $D-1$, we apply the Even Saturating Lemma (mapping $C^i$ to $\C$, $y_i$ to $\alpha$, $y'_i$ to $\alpha'$ and $x_{(i \bmod{k})+1}$ to $\gamma$) and obtain an additional short cycle $\C^1$ in $\G$ such that $x_{(i \bmod{k})+1}$ is a repeat  in $\C^1$ of a neighbor $\mu\not\in C^i$ of $y'_i$, and $\C^1\cap C^i=\emptyset$. Since $x_{(i \bmod{k})+1}$ is saturated, we have that necessarily $\C^1=C^{(i \bmod{k})+1}$, which, in turn, implies $\mu=x'_{(i \bmod{k})+1}$. In other words, it follows that $y'_i\sim x'_{(i \bmod{k})+1}\in E(\G)$.

If, instead, $I_i$ is a $(D-1)$-path then $I_{(i \bmod{k})+1}$ must be a path of length smaller than $D-1$, otherwise there would not exist a path $I_j$ ($1\le j\le k$) of length smaller than $D-1$, contrary to our assumptions. Therefore, we can apply the above reasoning and deduce that $x'_{(i \bmod k)+1}\sim y'_i \in E(\G)$.

This way we obtain a subgraph $\Upsilon= \bigcup_{i=1}^{k} \big( I'_i \cup y'_i\sim x'_{(i \bmod k)+1} \big)=x'_1I'_1y'_1x'_2I'_2y'_2\ldots x'_kI'_ky'_kx'_1$
intersecting $C^i$ at $I'_i$ for $1\le i\le k$ (see Fig.~\ref{fig:RepeatCycle}($b$), where part of the subgraph $\Upsilon$ is highlighted in bold).

We next show that $\Upsilon$ must be indeed a cycle.

{\bf Claim 1.} $\Upsilon$ is a $2D$-cycle.

{\bf Proof of Claim 1.} First note that $\Upsilon$ is connected and that $|\Upsilon|\le 2D$. By Proposition \ref{prop:Girth}, unless $\Upsilon$ is a $2D$-cycle, $\Upsilon$ contains no short cycle. If the neighbors of $C$ are pairwise disjoint then $\Upsilon$ is a $2D$-cycle. Suppose that some neighbors of $C$ are non-disjoint and that $\Upsilon$ is not a cycle, then $\Upsilon$ is a tree.

Let $z \in C^\ell$ be an arbitrary leaf in $\Upsilon$. If the repeat path $I'_\ell=x'_\ell-y'_\ell$ had length greater than 0, then $z$ would have at least two neighbors in $\Upsilon$. Therefore, $I_\ell=C \cap C^\ell$ contains exactly one vertex, and thus, $x_\ell = y_\ell$ and $z = x'_\ell = y'_\ell$.

Recall we do addition modulo $k$ on the subscripts of the vertices and the superscripts of the cycles.

\begin{figure}[!ht]
\begin{center}
\makebox[\textwidth][c]{\includegraphics[scale=1]{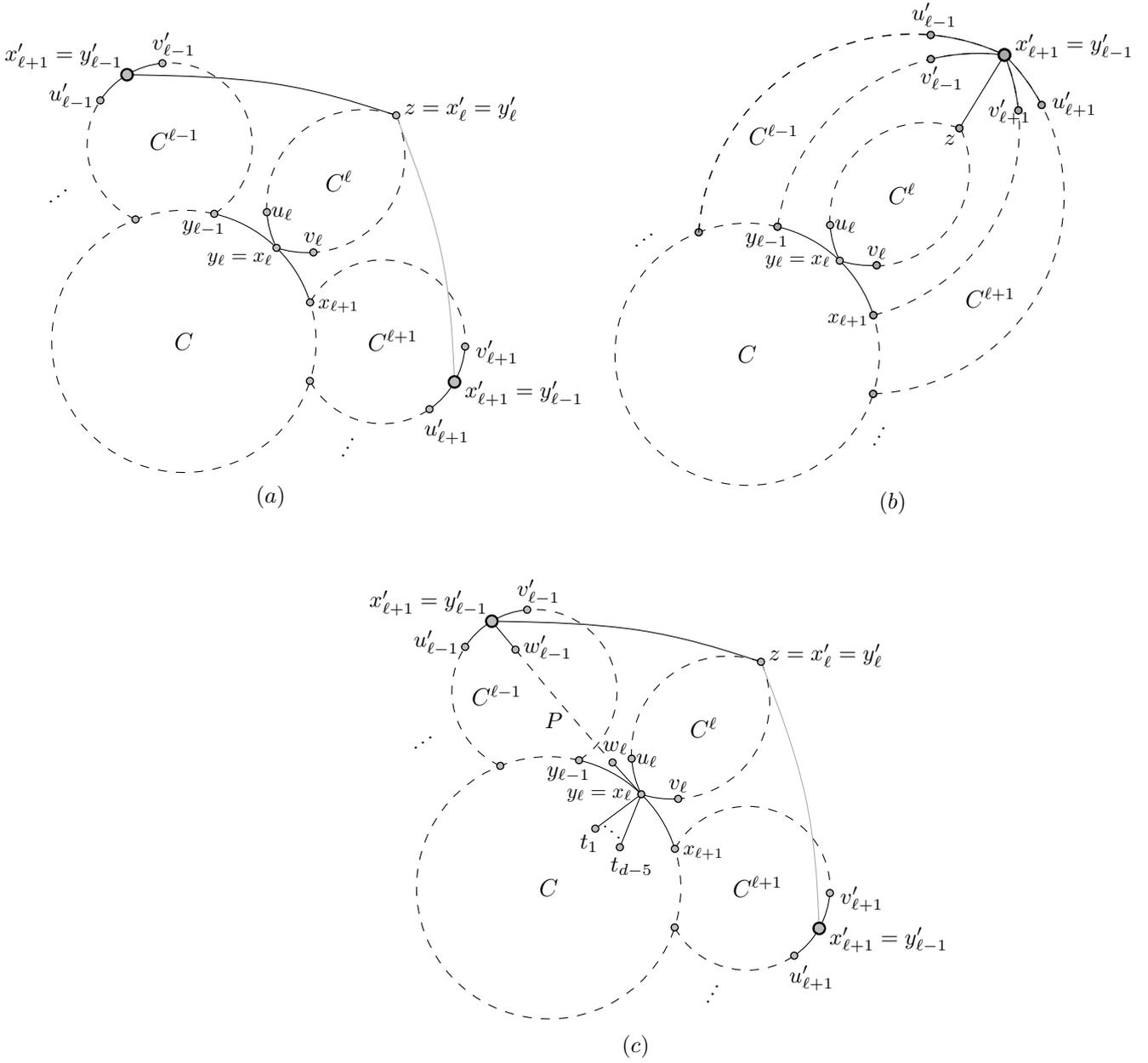}}
\caption{Auxiliary figure for Claim 1 of Lemma \ref{lem:RepeatCycle}.}
\label{fig:RepeatPatch}
\end{center}
\end{figure}

Since $x'_\ell \sim y'_{\ell-1}$ and $x'_\ell \sim x'_{\ell+1}$ are edges in $\Upsilon$, it holds that $y'_{\ell-1}$ and $x'_{\ell+1}$ denote the same vertex. Let $u'_{\ell-1}, v'_{\ell-1}$ be the neighbors of $y'_{\ell-1}$ in $C^{\ell-1}$; $u'_{\ell+1}, v'_{\ell+1}$ the neighbors of $x'_{\ell+1}$ in $C^{\ell+1}$; and $u_\ell, v_\ell$ the neighbors of $x_\ell$ in $C^\ell$. We have that $V(C^{\ell-1} \cap C^{\ell+1}) = \{y'_{\ell-1}\}$, otherwise there would be a third short cycle in $\G$ containing $x_\ell$. In particular, the vertices in $\{u'_{\ell-1}, v'_{\ell-1}, u'_{\ell+1}, v'_{\ell+1}, x'_\ell \}$ are pairwise distinct and $d\ge5$. See Fig.~\ref{fig:RepeatPatch} ($a$) and ($b$) for two drawings of this situation.

Now consider a path $P=x_\ell-y'_{\ell-1}$. Since $x_\ell$ cannot be contained in a further short cycle, we have that $P$ must be a $D$-path and go through a neighbor $w'_{\ell-1}$ of $y'_{\ell-1}$ not contained in $\{u'_{\ell-1}, v'_{\ell-1}, u'_{\ell+1}, v'_{\ell+1}, x'_\ell \}$, which implies $d \ge 6$. By similar arguments, we obtain that $P$ must go through a neighbor $w_\ell$ of $x_\ell$ not contained in $\{y_{\ell-1}, x_{\ell+1}, u_\ell, v_\ell \}$.

Finally, let $t_1,t_2,\ldots,t_{d-5}$ denote the vertices in $N(x_\ell)-\{y_{\ell-1}, x_{\ell+1}, u_\ell, v_\ell, w_\ell \}$; see Fig.~\ref{fig:RepeatPatch} ($c$). Consider a path $Q^i=t_i-y'_{\ell-1}$. Since $x_\ell$ cannot be contained in a further short cycle, $Q^i$ must be a $D$-path and go through a neighbor of $y'_{\ell-1}$ not contained in $\{u'_{\ell-1}, v'_{\ell-1}, u'_{\ell+1}, v'_{\ell+1}, x'_\ell, w'_{\ell-1} \}$. Therefore, we have that $d\ge7$ and, by the pigeonhole principle, that there are two paths $Q^r$ and $Q^s$ containing a common neighbor of $y'_{\ell-1}$. This way, $x_\ell$ would be contained in a third short cycle, a contradiction.

As a result, we conclude that the repeat graph $\Upsilon$ of $C$ is indeed a $2D$-cycle $C'$ as claimed. This completes the proof of Claim 1, and thus, of the lemma. \EndProof

We call the aforementioned cycle $C'$ the {\it repeat} of the cycle $C$ in $\G$, and denote it
by $\rep(C)$. Some simple consequences of the Repeat Cycle Lemma follow next.

\begin{corollary}[Repeat Cycle Uniqueness]
\label{cor:RepeatUniqueness}
If a $2D$-cycle $C$ has a repeat cycle $C'$ then $C'$ is unique.
\end{corollary}

\begin{corollary}[Repeat Cycle Symmetry]
\label{cor:RepeatSymmetry}
If $C'=\rep(C)$ then $C=\rep(C')$.
\end{corollary}

\begin{corollary}
\label{cor:RepeatPath}
Let $\G$ be a $(d,D,-2)$-graph with $d\ge 4$ and $D\ge 2$. Let $C,C^1$ be two $2D$-cycles in
$\G$ which intersect at a path $I$ of length smaller than $D-1$, and set $I'=\rep^{C^{1}}(I)$.
Then the repeat cycle of $C$ intersects $C^1$ at $I'$.
\end{corollary}

\begin{corollary}[Handy Corollary]
\label{cor:RepeatHandy}
Let $\G$ be a $(d,D,-2)$-graph with $d\ge 4$ and $D\ge 2$, $\C$ a $2D$-cycle in $\G$, and $x,x'$ repeat vertices in $\C$. Let $\C^1$ and $\C^2$ be $2D$-cycles other than $\C$ containing $x$ and $x'$, respectively. Suppose that $I=\C^1\cap \C$ is a path of length smaller than $D-1$. Then, setting $y=\rep^{\C^1}(x)$ and $y'=\rep^{\C^2}(x')$, we have that $y$ and $y'$ are repeat vertices in the repeat cycle of $\C$.
\end{corollary}

\pr We denote the $k$ neighbor cycles of $\C$ by $E^1,E^2,\ldots E^k$ and their respective
intersection paths with $\C$ by $I_1=x_1-y_1,I_2=x_2-y_2,\ldots,I_k=x_k-y_k$ in such a way that
$\C=x_1I_1y_1x_2I_2y_2\ldots x_kI_ky_kx_1$. For $1\le j\le k$, we also denote $I'_j=x'_j-y'_j$,
where $x'_j=\rep^{E^j}(x_j)$ and $y'_j=\rep^{E^j}(y_j)$.

Obviously, for some $r,s$ $(1\le r,s\le k)$ we have that $\C^1=E^r$, $\C^2=E^s$, $x\in I_r$,
$x'\in I_s$, $y\in I'_r$, and $y'\in I'_s$. We may assume $r<s$. By the Repeat Cycle Lemma, the
vertices $y$ and $y'$ belong to the repeat cycle $\C'$ of $\C$. Then the paths $xI_ry_rx_{r+1}I_{r+1}y_{r+1}\ldots x_{s-1}I_{s-1}y_{s-1}x_sI_sx'\subset \C$ and $yI'_ry'_rx'_{r+1}I'_{r+1}y'_{r+1}\ldots x'_{s-1}I'_{s-1}y'_{s-1}x'_sI'_sy'\subset \C'$ are both $D$-paths in $\G$, and the corollary follows.\EndProof


\section{Main Results}
\label{sec:MainResults}

\subsection{On the girth of $(d,D,-2)$-graphs}

\begin{proposition}
\label{prop:(2D-1)}
A $(d,D,-2)$-graph $\G$ with $d\ge 4$ and $D\ge 4$ does not contain $(2D-1)$-cycles.
\end{proposition}

\pr Suppose, by way of contradiction, that there is a $(2D-1)$-cycle $C$ in $\G$.

\begin{figure}[!ht]
\begin{center}
\makebox[\textwidth][c]{\includegraphics[scale=.9]{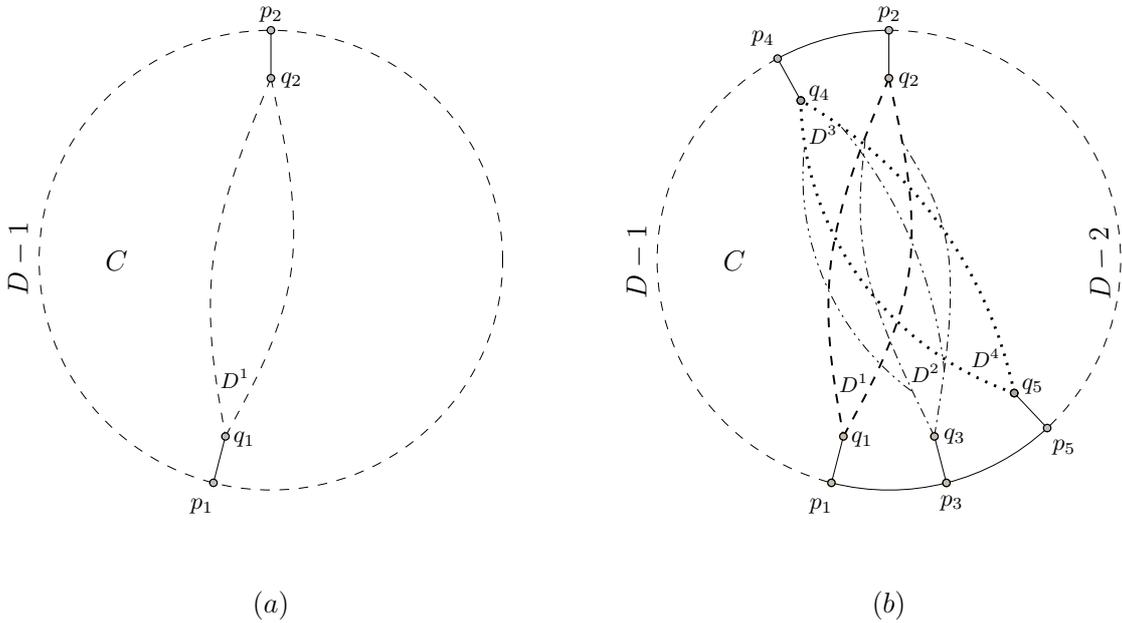}}
\caption{Auxiliary figure for Proposition \ref{prop:(2D-1)}}
\label{fig:(2D-1)}
\end{center}
\end{figure}

Let $p_1,p_2$ be two repeat vertices in $C$, and $q_1$ a neighbor of $p_1$ not contained in $C$. According to the Odd Saturating Lemma, there are both a neighbor $q_2$ of $p_2$ not contained in $C$ and a $2D$-cycle $D^1$, such that $q_1$ and $q_2$ are repeats in $D^1$ (see Fig.~\ref{fig:(2D-1)} $(a)$).

For $1\le i\le 3$, denote by $p_{i+2}$ the repeat of $p_{i+1}$ in $C$ other than $p_i$. We now
apply the Odd Saturating Lemma (mapping $C$ to $\C$, $p_2$ to $\alpha$, $p_{3}$ to $\alpha'_2$,
 $q_2$ to $\gamma$) and ascertain the existence of a $2D$-cycle $D^2$ and a neighbor
$q_{3}$ of $p_{3}$ not contained in $C$, such that $q_2$ and $q_{3}$ are repeats in $D^2$.
For $i=3,4$ by repeatedly applying the Odd Saturating Lemma (mapping $C$ to $\C$, $p_i$ to
$\alpha$, $p_{i+1}$ to $\alpha'_2$, $q_i$ to $\gamma$) we ensure the existence of a
$2D$-cycle $D^i$ and a neighbor $q_{i+1}$ of $p_{i+1}$ not contained in $C$, such that $q_i$ and $q_{i+1}$ are repeats in $D^i$. See Fig.~\ref{fig:(2D-1)} $(b)$.

Note that $D^1\cap D^2$ is a path of length at most $2<D-1$; otherwise for some vertex $t\in D^1\cap D^2$ the cycle $tD^1q_1p_1p_3q_3D^2t$ would have length at most $2D-1$, a contradiction. Similarly, $D^2\cap D^3$ and $D^3\cap D^4$ are paths of length at most $2$.

We now apply the Handy Corollary. By mapping the cycle $D^2$ to $\C$, the vertex $q_2$ to $x$, the vertex $q_3$ to $x'$, the cycle $D^1$ to $\C^1$, the cycle $D^3$ to $\C^2$, the vertex $q_1$ to $y$ and the vertex $q_4$ to $y'$, we obtain that $q_1$ and $q_4$ are repeat vertices in the repeat cycle of $D^2$. Therefore, since $q_4\in D^4$,  it follows that $D^2$ and $D^4$ are repeat cycles and $q_1=q_5$. As a consequence, there is in $\G$ a cycle $q_1p_1p_3p_5q_5$ of length $4<2D-1$, a contradiction. \EndProof

From Propositions \ref{prop:Girth} and \ref{prop:(2D-1)}, it  follows immediately that

\begin{theorem}
\label{theo:Girth2D}
The girth of a $(d,D,-2)$-graph $\G$ with $d\ge 4$ and $D\ge 4$ is $2D$.
\end{theorem}

\subsection{Non-existence of subgraphs isomorphic to $\Theta_D$}

\begin{proposition}
\label{prop:Theta}
A $(d,D,-2)$-graph $\G$ with $d\ge 4$ and $D\ge 4$ does not contain a subgraph isomorphic to $\Theta_D$.
\end{proposition}

\pr In this proof our reasoning resembles that of Proposition \ref{prop:(2D-1)}, and especially,
 of \cite[Proposition 5.1 ]{FP10}.

\begin{figure}[!ht]
\begin{center}
\makebox[\textwidth][c]{\includegraphics[scale=.9]{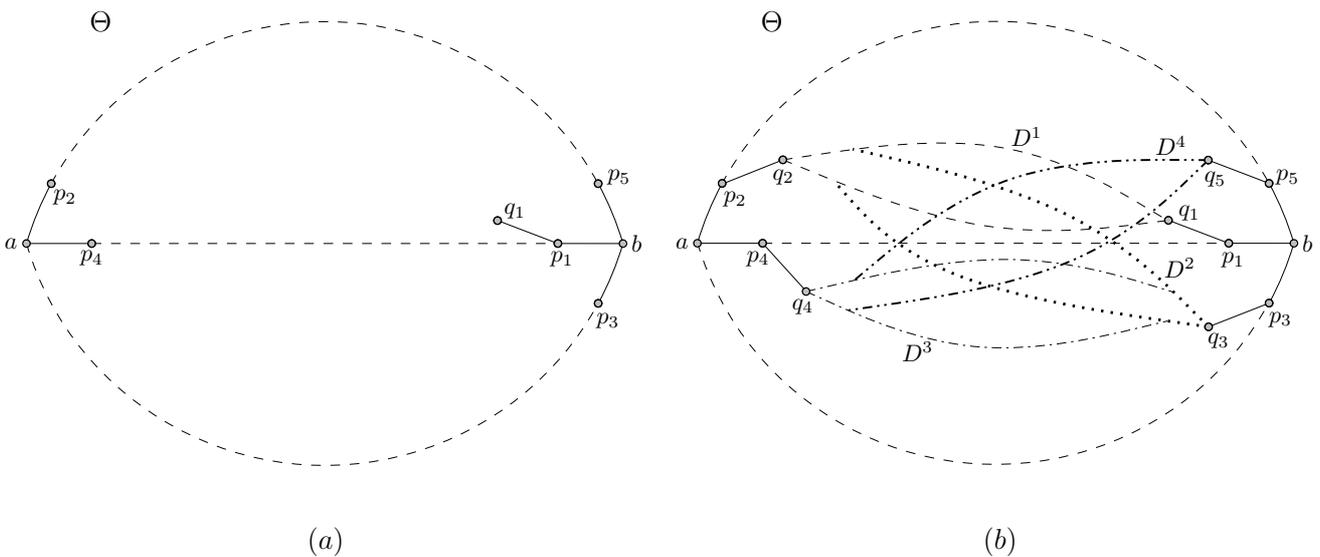}}
\caption{Auxiliary figure for Proposition \ref{prop:Theta}}
\label{fig:Theta}
\end{center}
\end{figure}

Suppose that $\G$ contains a subgraph $\Theta$ isomorphic to $\Theta_D$, with branch vertices $a$ and $b$. Let $p_1,p_2,p_3,p_4$ and $p_5$ be as in Fig.~\ref{fig:Theta} ($a$), and let $q_1$ be one of the neighbors of $p_1$ not contained in $\Theta$.

Since all vertices of $\Theta$ are saturated, there cannot be a short cycle in $\G$ containing any of the incident edges of $p_1,p_2,p_3,p_4$ or $p_5$  which are not contained in $\Theta$. According to this and by applying the Even Saturating Lemma, there is an additional $2D$-cycle $D^1$ in $\G$ such that $q_1$ and one of the neighbors of $p_2$ not contained in $\Theta$ (say $q_2$) are repeats in $D^1$. Also, it follows that $D^1\cap \Theta=\emptyset$. Analogously, by repeatedly applying the Even Saturating Lemma, for $2\le i\le 4$ we obtain that there is an additional $2D$-cycle $D^i$ such that $q_i$ and one of the neighbors of $p_{i+1}$ not contained in $\Theta$ (say $q_{i+1}$) are repeats in $D^i$. Also, we have that $D^i\cap \Theta=\emptyset$ (see Fig.~\ref{fig:Theta} ($b$)).

Note that $D^1\cap D^2$ is a path of length at most $2<D-1$; otherwise for some vertex $t\in D^1\cap D^2$ there would be a cycle $tD^1q_1p_1bp_3q_3D^2t$ of length at most $2D$ to which the vertex $b$ would belong, a contradiction. For similar reasons, the intersection paths $D^2\cap D^3$ and $D^3\cap D^4$ have both length at most 2.

We now apply the Handy Corollary. By mapping the cycle $D^2$ to $\C$, the vertex $q_2$ to $x$,
the vertex $q_3$ to $x'$, the cycle $D^1$ to $\C^1$, the cycle $D^3$ to $\C^2$, the vertex
$q_1$ to $y$ and the vertex $q_4$ to $y'$, we obtain that $q_1$ and $q_4$ are repeat vertices
in the repeat cycle of $D^2$. Therefore, since $q_4\in D^4$,  it follows that $D^2$ and $D^4$
are repeat cycles and $q_1=q_5$; but then there is a cycle $q_1p_1bp_5q_5$ in $\G$ of length $4<2D$, a contradiction.\EndProof

\begin{corollary}
\label{cor:Type(iii)}
Every vertex in a $(d,D,-2)$-graph $\G$ with $d\ge 4$ and $D\ge 4$ is of Type $(iii)$.
\end{corollary}

\subsection{Non-existence results on $(d,D,-2)$-graphs}

In view of Corollary \ref{cor:Type(iii)}, the following corollary, which was obtained in \cite{DP10}, follows immediately.

\begin{corollary}[Corollary 2.3 from \cite{DP10}]\label{congruences}
The feasible values of $d$ for $(d,D,-2)$-graphs are restricted
according to the following conditions.

\begin{itemize}
\item[] When $D$ is even, $d$ is odd.
\item[] When $D$ is a power of an odd prime, $d-1$ is a multiple of $D$.
\item[] When $D\ge4$ is a power of 2,  $d-1$ is a multiple of $D/2$.
\end{itemize}

\end{corollary}

\begin{proposition}
\label{prop:Number2DCycles}
The number $N_{2D}$ of $2D$-cycles in a $(d,D,-2)$-graph $\G$ with $d\ge 4$ and $D\ge 4$ is given by the expression $N_{2D}=\frac{n}{D}=\frac{d\big(1+(d-1)+\ldots+{(d-1)}^{D-1}\big)-1}{D}$, where $n$ is the order of $\G$.
\end{proposition}

\pr According to Proposition \ref{prop:Girth} and Corollary \ref{cor:Type(iii)}, every vertex of $\G$ is contained in exactly two $2D$-cycles. We
then count the number $N_{2D}$ of $2D$-cycles of $\G$ . Since the order of $\G$ is
$n=1+d\big(1+(d-1)+\ldots+{(d-1)}^{D-1}\big)-2$, we have that

\begin{center}
$N_{2D}=\frac{2\times\Big(1+d\big(1+(d-1)+\ldots+{(d-1)}^{D-1}\big)-2\Big)}{2D}=\frac{d\big(1+(d-1)+\ldots+{(d-1)}^{D-1}\big)-1}{D}$,
\end{center}

and the proposition follows. \EndProof

\begin{lemma}
\label{lem:IntersectionLength}
Every two non-disjoint $2D$-cycles in a $(d,D,-2)$-graph $\G$ with $d\ge 4$ and $D\ge 4$ intersect at a path of length smaller than $D-1$.
\end{lemma}

\pr We follow a strategy very similar to the one used in the proof of \cite[Lemma 5.1]{FP10}.

Since $\G$ does not contain a graph isomorphic to $\Theta_{D}$, it is only necessary to prove here that any two non-disjoint $2D$-cycles in $\G$ cannot intersect at a path of length $D-1$. Suppose, by way of contradiction, that there are two $2D$-cycles $C^1$ and $C^2$ in $\G$ intersecting at a path $I_1$ of length $D-1$.

Let $v$ be an arbitrary vertex on $I_1$, and $v'=\rep^{C^2}(v)$. Let $C^3$ be the other $2D$-cycle containing $v'$, and $I_2=C^2\cap C^3$. If $I_2$ were a path of length smaller than $D-1$
 then, by Corollary \ref{cor:RepeatPath}, the repeat cycle of $C^3$ would intersect $C^2$ at a proper subpath of $I_1$ containing $v$. This is a clear contradiction to the fact that $v$ is already saturated. Consequently, $I_2$ must be a $(D-1)$-path and $C^2$ is intersected by exactly two $2D$-cycles, namely $C^1$ and $C^3$, at two independent $(D-1)$-paths.

\begin{figure}[!ht]
\begin{center}
\makebox[\textwidth][c]{\includegraphics[scale=.9]{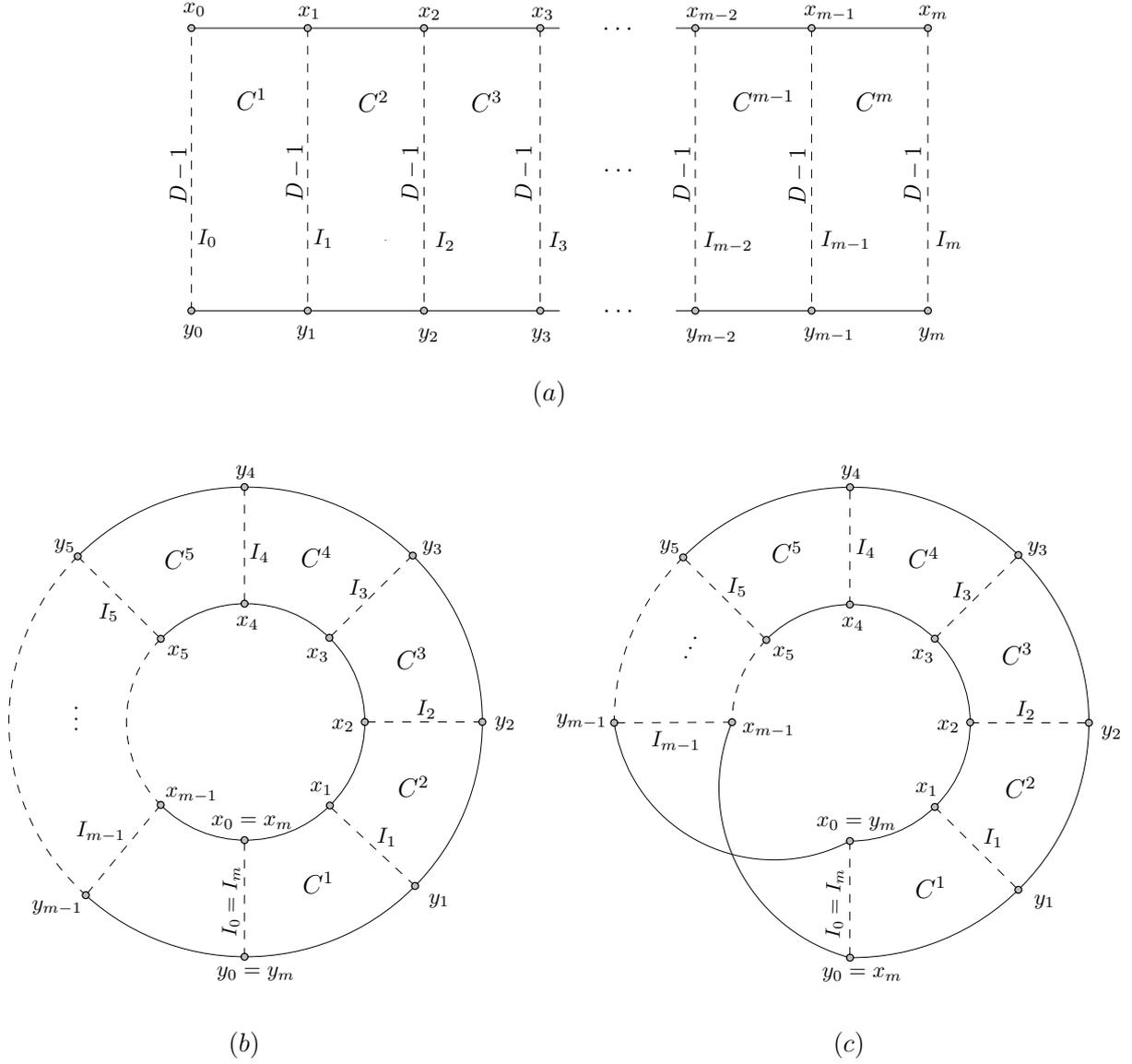}}
\caption{Auxiliary figure for Lemma \ref{lem:IntersectionLength}}
\label{fig:IntersectionLength1}
\end{center}
\end{figure}

By repeatedly applying this reasoning and considering that $\G$ is finite, we obtain a maximal length sequence $C^1,C^2,C^3,\ldots,C^{m}$ of pairwise disctinct $2D$-cycles in $\G$ such that $C^i$ intersects $C^{i+1}$ at a path $I_i$ of length $D-1$ ($1\le i\le m-1$). Furthermore, it follows that $C^j\cap C^k=\emptyset$ for any $j,k\in\{1,\ldots,m\}$ such that $2\le|i-j|\le m-2$. Let us denote the paths $I_1=x_1-y_1,\ldots,$ $I_{m-1}=x_{m-1}-y_{m-1}$ in
such a way that, for $1\le i\le m-2$, $x_i\sim x_{i+1}$ and $y_{i}\sim y_{i+1}$ are edges in $\G$. Also, let $x_0\in N(x_{1})\cap (C^1-I_1)$, $y_0\in N(y_{1})\cap (C^1-I_1)$, $x_m\in N(x_{m-1})\cap (C^m-I_{m-1})$, and $y_m\in N(y_{m-1})\cap (C^m-I_{m-1})$. Figure~\ref{fig:IntersectionLength1} $(a)$ shows this configuration. Set $I_0=x_0-y_0$ and $I_m=x_m-y_m$. Since the sequence $C^1,C^2,C^3,\ldots,C^{m}$ is maximal and all the vertices in $I_1,\ldots,I_{m-1}$ are saturated, it follows that $I_0=I_m$, and we have either $x_0=x_m$ and $y_0=y_m$ (as in Fig.~\ref{fig:IntersectionLength1} $(b)$), or $x_0=y_m$ and $y_0=x_m$ (as in Fig.~\ref{fig:IntersectionLength1} $(c)$).

If $x_0=x_m$ and $y_0=y_m$ then $m>2D$; otherwise the cycle $x_1x_2\ldots x_mx_1$ would have length at most $2D$, contradicting the saturation of $x_1$. If, conversely, $x_0=y_m$ and $y_0=x_m$ then $m>D$; otherwise the cycle $x_1x_2\ldots x_my_1y_2\ldots y_mx_1$ containing $x_1$ would have length at most $2D$, a contradiction. For our purposes, it is enough to state $m>D\ge4$ in any case.

We now proceed with the second part of the proof.

\begin{figure}[!ht]
\begin{center}
\makebox[\textwidth][c]{\includegraphics[scale=.9]{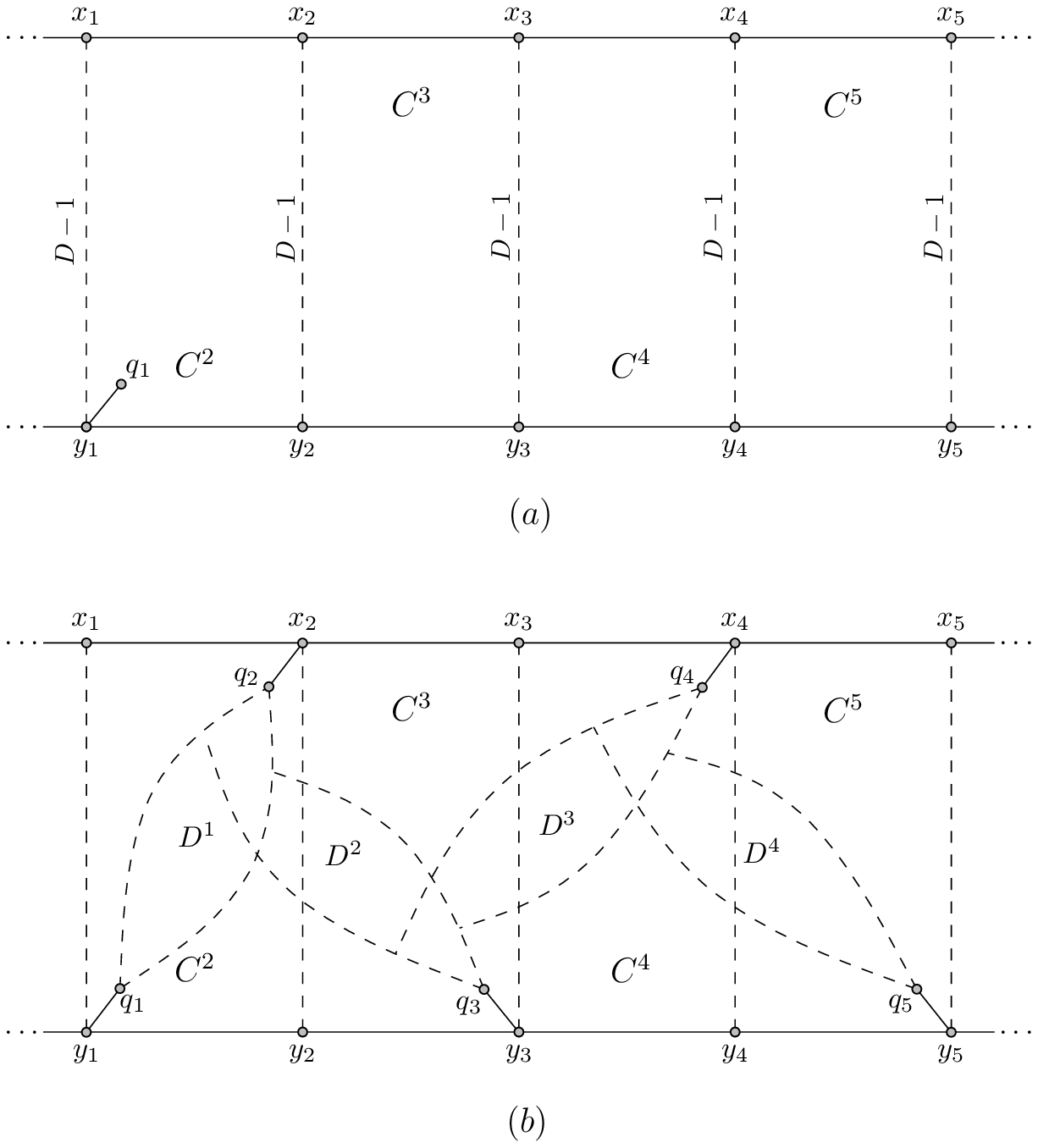}}
\caption{Auxiliary figure for Lemma \ref{lem:IntersectionLength}}
\label{fig:IntersectionLength2}
\end{center}
\end{figure}

Let $\Phi=\cup_{i=1}^m C^i$, and $q_1$ a neighbor of $y_1$ not contained in $\Phi$ (see Fig.~\ref{fig:IntersectionLength2} ($a$)).

Since $y_1$ is saturated, the edge $q_1\sim y_1$ cannot be contained in a further short cycle. We apply the Even Saturating Lemma (by mapping $C^{2}$ to $\C$, $y_1$ to $\alpha$, $x_2$ to $\alpha'$, and $q_1$ to $\gamma$), and obtain in $\G$ an additional $2D$-cycle $D^1$ such that $q_1$ and one of the neighbors of $x_2$ not contained in $\Phi$ (say $q_2$) are repeats in $D^1$, and $D^1\cap C^2=\emptyset$.
 Analogously, there exists an additional $2D$-cycle $D^2$ such that $q_2$ and a neighbor of $y_3$ not contained in $\Phi$ (say $q_3$) are repeats in $D^2$, and $D^2\cap C^3=\emptyset$; an additional $2D$-cycle $D^3$ such that $q_3$ and a neighbor of $x_4$ not contained in $\Phi$ (say $q_4$) are repeats in $D^3$, and $D^3\cap C^4=\emptyset$; and an additional $2D$-cycle $D^4$ such that $q_4$ and a neighbor of $y_5$ not contained in $\Phi$ (say $q_5$) are repeats in $D^4$, and $D^4\cap C^5=\emptyset$. See Fig.~\ref{fig:IntersectionLength2} ($b$).

Note that $D^1\cap D^2$ cannot be a $(D-1)$-path; otherwise for some vertex $t\in D^1\cap D^2$ there would be a cycle $tD^1q_1y_1y_2y_3q_3D^2t$ of length at most $4+D-2+D-2$ (since $D-1\ge 3$), a contradiction to the fact that $y_1$ is saturated and $\g(\G)=2D$. Analogously, $D^i\cap D^{i+1}$ cannot be a $(D-1)$-path for $i=2,3$.

We now apply the Handy Corollary as in the proofs of the previous theorems. By mapping  the cycles $D^2$ to $\C$, $D^1$ to $\C^1$ and $D^3$ to $\C^2$, and the vertices $q_2$ to $x$, $q_3$ to $x'$, $q_1$ to $y$, and $q_4$ to $y'$, it follows that the vertices $q_1$ and $q_4$ are repeat vertices in the repeat cycle of $D^2$. Since $q_4\in D^4$, we have that $D^2$ and $D^4$ are repeat cycles and that $q_5=q_1$. This way, we obtain a cycle $q_1y_1y_2y_3y_4y_5q_5$ in $\G$ of length $6<2D$, a contradiction.

This completes the proof of the lemma.\EndProof

We are now in a position to prove our second main result.


\begin{theorem}
\label{theo:Main1}
There are no $(d,D,-2)$-graphs with even $d\ge 4$ and $D\ge 4$.
\end{theorem}

\pr Suppose there is a $(d,D,-2)$-graph $\G$ with even $d\ge 4$ and $D\ge 4$.

According to Lemma \ref{lem:IntersectionLength}, any two non-disjoint $2D$-cycles in $\G$ intersect at a path of length smaller than $D-1$, which means that every $2D$-cycle $C$ in $\G$ has a repeat cycle $C'$ (by the Repeat Cycle Lemma). Because of the uniqueness and symmetry of repeat cycles, the number $N_{2D}$ of $2D$-cycles in $\G$ must be even.

However, since $d$ is even, the number $N_{2D}=\frac{d\big(1+(d-1)+\ldots+{(d-1)}^{D-1}\big)-1}{D}$ of $2D$-cycles in $\G$ is odd, a contradiction.\EndProof

Note that Theorem \ref{theo:Main1} contains, as a special case, the result of the non-existence
of $(4,D,-2)$-graphs for $D\ge4$, which was claimed prematurely in \cite{MS05}.

From Proposition \ref{prop:Number2DCycles} we easily derive the following results:

\begin{theorem}
\label{theo:Main2}
There are no $(d,D,-2)$-graphs with odd $d\ge 5$, $D\ge 4$ and order $n=d\big(1+(d-1)+\ldots+{(d-1)}^{D-1}\big)-1\not\equiv0\pmod D$.
\end{theorem}

\begin{corollary}
\label{theo:Main3}
There are no $(d,D,-2)$-graphs with odd $d\ge 5$ and $D\ge 5$ such that $d\equiv0,2\pmod D$.
\end{corollary}

Furthermore, for a particular value of $D\ge5$ it is possible to rule out the existence of $(d,D,-2)$-graphs with odd $d\ge5$ for many other values of $d$, by considering the set of all possible residues of $d$ in the division by $D$. If, for some $r\in\{0,1,\ldots,D-1\}$, we have $d\equiv r\pmod D$ implies $d\big(1+(d-1)+\ldots+{(d-1)}^{D-1}\big)-1\not\equiv0\pmod D$, then there are no $(d,D,-2)$-graphs with odd $d\ge5$ such that $d\equiv r\pmod D$.

Accordingly, the following table shows all values of $4\le D\le 16$ and odd $d\ge 5$ for which a $(d,D,-2)$-graph  might still exist.

\begin{center}
\begin{tabular}{|l|l|}
\hline
{\bf $D$}	&	{\bf $d$}				\\
\hline
$4$		&	$d\equiv1,3\pmod4$		\\
$5$		&	$d\equiv1\pmod{10}$		\\
$6$		&	$d\equiv1\pmod6$			\\
$7$		&	$d\equiv1\pmod{14}$		\\
$8$		&	$d\equiv1,5\pmod8$		\\
$9$		&	$d\equiv1\pmod{18}$		\\
$10$		&	$d\equiv1,9\pmod{10}$		\\
$11$		&	$d\equiv1\pmod{22}$		\\
$12$		&	$d\equiv1,7\pmod{12}$		\\
$13$		&	$d\equiv1\pmod{26}$		\\
$14$		&	$d\equiv1,13\pmod{14}$		\\
$15$		&	$d\equiv1,13\pmod{30}$		\\
$16$		&	$d\equiv1,9\pmod{16}$		\\
\hline
\end{tabular}
\end{center}

\section{Non-existence of $(4,3,-2)$-graphs}

In this section we prove the non-existence of $(4,3,-2)$-graphs (see Theorem
\ref{theo:(4,3,-2)}), which will allow us to provide the full catalogue of $(4,D,-2)$-graphs with $D\ge2$.

Proposition \ref{prop:Theta} asserts the non-existence of a subgraph isomorphic to
$\Theta_D$ in a $(d,D,-2)$-graph $\G$ with $d\ge 4$ and $D\ge 4$. We next give an alternative proof for $d=4$ that covers also the case $D=3$.

\begin{proposition}
\label{prop:Theta_D} A $(4,D,-2)$-graph $\G$ with $D\ge 3$ does not contain a
subgraph isomorphic to $\Theta_{D}$.
\end{proposition}

\begin{figure}[phtb]
\begin{center}
\makebox[\textwidth][c]{\includegraphics[scale=1]{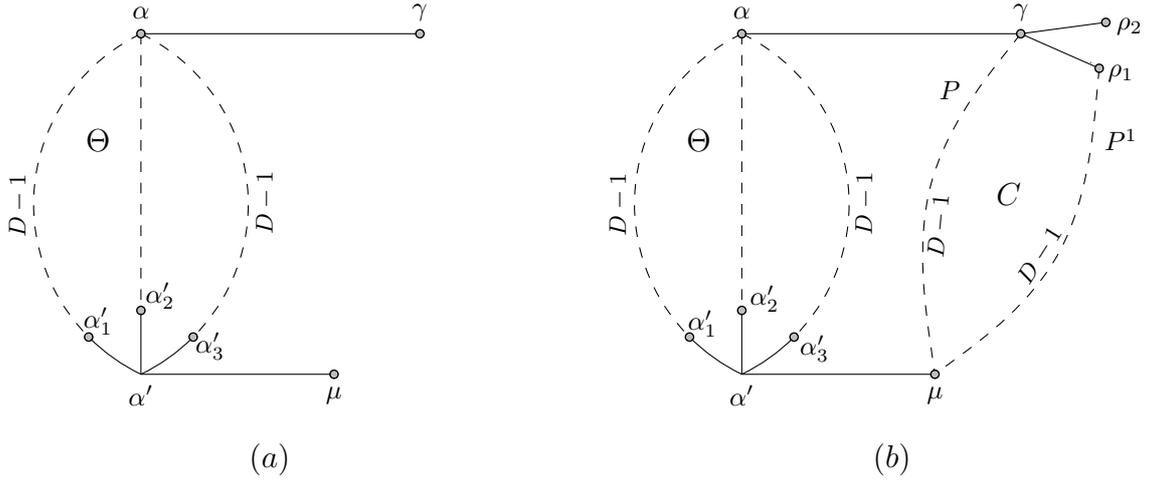}}
\caption{Auxiliary figure for Proposition \ref{prop:Theta_D}.}
\label{fig:Theta_D}
\end{center}
\end{figure}

\pr Suppose that $\G$ contains a subgraph $\Theta$ isomorphic to $\Theta_D$, where $\alpha$ and $\alpha'$ are its branch vertices. Let $\alpha'_1$, $\alpha'_2$, $\alpha'_3$, $\gamma$ and $\mu$ be as in Figure~\ref{fig:Theta_D} ($a$).

First consider a path $P=\gamma-\alpha'$. As $\alpha$ cannot belong to any further short cycle, $P$ must go through $\mu$ and be a $D$-path. Let $\rho_1$ and $\rho_2$ be the neighbors of $\gamma$ other than $\alpha$ and not contained in $P$. Consider a path $P^1=\rho_1-\alpha'$. As $\alpha$ is saturated, $P^1$ cannot go through $\alpha'_1$, $\alpha'_2$ or $\alpha'_3$, so it must go through $\mu$ and be a $D$-path. This way, $\gamma$ is contained in a $(2D-1)$-cycle $C=\gamma P\mu P^1\rho_1\gamma$, and $\gamma$ becomes saturated. Analogously, a path $P^2=\rho_2-\alpha'$ must go through $\mu$, causing the formation of another short cycle containing $\mu$, a contradiction to Proposition \ref{prop:Girth} $(ii)$. See Figure~\ref{fig:Theta_D} ($b$). \EndProof

Next we prove that the girth of a $(4,3,-2)$-graph must be 6 by ruling out the existence of $5$-cycles.

\begin{proposition}
\label{prop:5cycle} A $(4,3,-2)$-graph $\G$  has girth $6$.
\end{proposition}

\pr We proceed by contradiction, supposing there is a 5-cycle $C$ in $\G$. In view of Proposition \ref{prop:Girth}, the graph $\G$ contains the subgraph $G$ of Fig.~\ref{fig:(4,3,-2)}, where $T_i$ denotes the enclosed set of 6 vertices at distance 2 from $x_i$, for $1\le i\le 5$.

\begin{figure}[phtb]
\begin{center}
\makebox[\textwidth][c]{\includegraphics[scale=1]{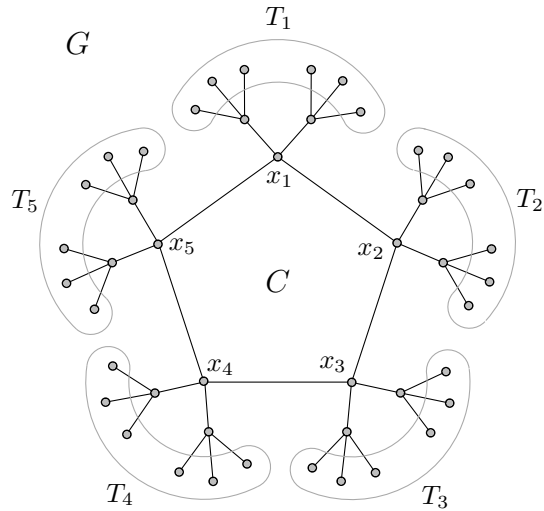}}
\caption{Auxiliary figure for Proposition \ref{prop:5cycle}.}
\label{fig:(4,3,-2)}
\end{center}
\end{figure}

Since $|\G|=51$ and $|G|=45$, there is a set $X\subset V(\G)$ such that $|X|=6$ and
$X\cap V(G)=\emptyset$. Any vertex $x\in X$ must be adjacent to a vertex in $T_i$, for
$1\le i\le 5$, in order to reach $x_i$ in at most 3 steps. However, this is clearly impossible since $\G$ has degree $4$. \EndProof

In view of Propositions \ref{prop:Girth}, \ref{prop:Theta_D}, \ref{prop:(2D-1)} and
\ref{prop:5cycle}, it follows that every vertex in a $(4,D,-2)$-graph $\G$ with $D\ge3$ is
contained in exactly two short cycles, namely, two $2D$-cycles.

\begin{proposition}
\label{prop:2D+1Cycles} The number $N_{2D+1}$ of $(2D+1)$-cycles in a $(4,D,-2)$-graph $\G$ with
$D\ge3$ is given by $N_{2D+1} = \frac{2\times3^{D}(2\times3^{D}-3)}{2D+1}$.
\end{proposition}

\pr The number of $(2D+1)$-cycles in $\G$ is closely
related to the number of edges involving only vertices at distance $D$ from any vertex $x$ in  $\G$. The number of vertices
at level $D$ is $4\times3^{D-1}-2$, and the number of elements in the set $F$ of edges involving only vertices at distance $D$ from $x$ is $$|F|=\frac{2\times2+3(4\times3^{D-1}-4)}{2}=2\times3^{D}-4,$$ since $x$ is contained in exactly two $2D$-cycles $C^1$ and $C^2$.

Denote by $y_1$ and $y_2$ the vertices at distance $D$ from $x$ on $C^1$ and $C^2$, respectively. Before proceeding to count, we prove that $y_1\sim y_2\not\in E(\G)$.

{\bf Claim 1.} $y_1\sim y_2\not\in E(\G)$.

\begin{figure}[phtb]
\begin{center}
\makebox[\textwidth][c]{\includegraphics[scale=1]{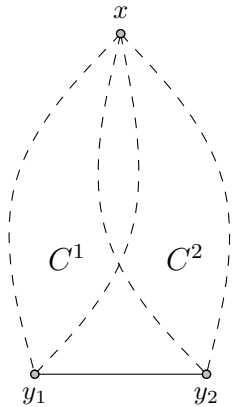}}
\caption{Auxiliary figure for Proposition \ref{prop:2D+1Cycles}.}
\label{fig:2D+1Cycles}
\end{center}
\end{figure}

{\bf Proof of Claim 1.} Suppose, by way of contradiction, that $y_1\sim y_2\in E(\G)$. Since $\g(\G)=2D$, it holds that $V(C^1\cap C^2)=\{x\}$; see Fig.~\ref{fig:2D+1Cycles}.  By Corollary \ref{cor:RepeatPath}, the repeat cycle $C'$ of $C^1$ intersects $C^2$ exactly at $y_2$; consequently, $C'$ contains the edge $y_1\sim y_2$. However, this contradicts the fact that $C^1$ and its repeat cycle $C'$ must be disjoint cycles.\EndProof

Accordingly, we partition the set $F$ into $F_1$, $F_2$ and $F_3$, where $F_1$ and $F_2$ are the sets of edges in $F$ adjacent to the vertices $y_1$ and $y_2$, respectively, and $F_3$ contains the remaining edges in $F$.

Each edge in $F_1$ or $F_2$ determines two $(2D+1)$-cycles containing $x$, while each edge from $F_3$ determines only one $(2D+1)$-cycle containing $x$. Therefore, given that $|F_1|=|F_2|=2$, we have that the number of $(2D+1)$-cycles passing through the vertex $x$ is

$$2|F_1|+2|F_2|+|F_3|=4+4+2\times3^{D}-8=2\times3^{D}.$$

Thus, the total number of $(2D+1)$-cycles in $\G$ is given by the expression

$$N_{2D+1}=\frac{2\times3^{D}(2\times3^{D}-3)}{2D+1},$$

and the proposition follows.\EndProof

Now we can readily prove Theorem \ref{theo:(4,3,-2)}.

\begin{theorem} There is no $(4,3,-2)$-graph. \label{theo:(4,3,-2)}
\end{theorem}

\pr By Proposition \ref{prop:2D+1Cycles}, the number of $7$-cycles in a $(4,3,-2)$-graph is
$2\times3^{3}(2\times3^{3}-3)/7=2754/7$, which is a contradiction. \EndProof


Theorems  \ref{theo:Main1} and \ref{theo:(4,3,-2)}  tell us that the $(4,2,-2)$-graph of
Fig.~\ref{fig:(Delta,D,-2)-graphs} ($d$) is the only  $(4,D,-2)$-graph for $D\ge2$. Thus, we
have successfully completed the census of all $(4,D,-2)$-graphs.

\section{Conclusions}
\label{sec:Conclusions}

In this paper, by exploiting the idea of extending the concept of repeats to paths and cycles,
put forward in \cite{FP10}, we obtained the results summarized below.

First, we proved that the girth of a $(d,D,-2)$-graph with $d\ge 4$ and $D\ge 4$ is $2D$. By obtaining necessary conditions for the existence of $(d,D,-2)$-graphs with $d\ge4$ and $D\ge4$, we proved the non-existence of $(d,D,-2)$-graphs with even $d\ge 4$ and $D\ge 4$. This outcome, together with a non-existence proof of $(4,3,-2)$-graphs, completed the catalogue of $(4,D,-\epsilon)$-graphs with $D\ge2$ and $0\le \epsilon\le 2$.

 \begin{description}
 \item[Catalogue of $(4,D,0)$-graphs with $D\ge2$.] There is no Moore graph of degree 4 and diameter $D\ge2$.
\item[Catalogue of $(4,D,-1)$-graphs with $D\ge2$.] There is no $(4,D,-1)$-graph for $D\ge2$.
\item[Catalogue of $(4,D,-2)$-graphs with $D\ge2$.] There is a unique $(4,2,-2)$-graph, shown in Fig.~\ref{fig:(Delta,D,-2)-graphs} ($d$).
\end{description}

We proved the non-existence of $(d,D,-2)$-graphs with odd $d\ge5$ and $D\ge 5$ such that
$d\equiv0,2\pmod{D}$. Furthermore, our new necessary conditions  allow us also to rule out the existence
of graphs of defect 2 for many other values of $d$ and $D$ using a simple approach.

\subsection{Remarks on the upper bound for $\N(\D,D)$}

Our results improve the upper bound on $N(\D,D)$ for many combinations of $\D$ and $D$.

\begin{proposition}
For even $\D\ge 4$ and $D\ge4$, $\N(\D,D)\le \M(\D,D)-3$.
\end{proposition}

In the particular case of $\D=4$, we have that $\N(4,2)=\M(4,2)-2$ and  $\N(4,D)\le \M(4,D)-3$ for $D\ge 3$.

According to Proposition \ref{prop:regularity}, a $(\D,D,-3)$-graph $\G$ must be regular; consequently, $(\D,D,-3)$-graphs with odd $\D\ge5$ and $D\ge4$ do not exist.

\begin{proposition}
For odd $\D\ge 5$ and $D\ge4$ such that $\D\big(1+(\D-1)+\ldots+{(\D-1)}^{D-1}\big)-1\not\equiv0\pmod{D}$, $\N(\D,D)\le \M(\D,D)-4$.
\end{proposition}

\begin{corollary}
For odd $\D\ge 5$ and $D\ge5$ such that $\D\equiv0,2\pmod{D}$, $\N(\D,D)\le \M(\D,D)-4$.
\end{corollary}

Finally, we feel that the following conjectures also hold.

\begin{conjecture}
There are no $(\D,D,-2)$-graph with $\D\ge 4$ and $D\ge 4$.
\end{conjecture}

\begin{conjecture}
For odd $\D\ge 5$ and $D\ge 4$, $\N(\D,D)\le \M(\D,D)-4$.
\end{conjecture}


\def\cprime{$'$}
\providecommand{\bysame}{\leavevmode\hbox to3em{\hrulefill}\thinspace}
\providecommand{\MR}{\relax\ifhmode\unskip\space\fi MR }
\providecommand{\MRhref}[2]{%
  \href{http://www.ams.org/mathscinet-getitem?mr=#1}{#2}
}
\providecommand{\href}[2]{#2}


\begin{thebibliography}{10}

\bibitem{BI73}
E.~Bannai and T.~Ito, \emph{On finite {M}oore graphs}, Journal of the Faculty
  of Science. University of Tokyo. Section IA. Mathematics \textbf{20} (1973),
  191--208.

\bibitem{BI81}
\bysame, \emph{Regular graphs with excess one}, Discrete Mathematics
  \textbf{37} (1981), no.~2-3, 147--158,
  \href{http://dx.doi.org/10.1016/0012-365X(81)90215-6}{doi:10.1016/0012-365X(%
81)90215-6}.

\bibitem{BA88}
H.~J. Broersma and A.~A. Jagers, \emph{The unique {$4$}-regular graphs on
  {$14$} and {$15$} vertices with diameter {$2$}}, Ars Combinatoria \textbf{25}
  (1988), no.~C, 55--62, Eleventh British Combinatorial Conference (London,
  1987).

\bibitem{Bus00}
D.~Buset, \emph{Maximal cubic graphs with diameter 4}, Discrete Applied
  Mathematics \textbf{101} (2000), no.~1-3, 53--61,
  \href{http://dx.doi.org/10.1016/S0166-218X(99)00204-8}{doi:10.1016/S0166-218%
X(99)00204-8}.

\bibitem{CG08}
J.~Conde and J.~Gimbert, \emph{On the existence of graphs of diameter two and
  defect two}, Discrete Mathematics \textbf{309} (2009), no.~10, 3166--3172,
  \href{http://dx.doi.org/10.1016/j.disc.2008.09.017}{doi:10.1016/j.disc.2008.%
09.017}.

\bibitem{Da73}
R.~M. Damerell, \emph{On {M}oore graphs}, Mathematical Proceedings of the
  Cambridge Philosophical Society \textbf{74} (1973), 227--236.

\bibitem{DP10}
C.~Delorme and G.~Pineda-Villavicencio, \emph{On graphs with cyclic defect or
  excess},  (2010), submitted.

\bibitem{Die05}
R.~Diestel, \emph{{Graph Theory}}, 3rd. ed., Graduate Texts in Mathematics,
  vol. 173, Springer-Verlag, Berlin, 2005.

\bibitem{El64}
B.~Elspas, \emph{Topological constraints on interconnection-limited logic},
  Proceedings of the Fifth Annual Symposium on Switching Circuit Theory and
  Logical Design (Princeton, NJ), IEEE Computer Society, 1964,
  \href{http://doi.ieeecomputersociety.org/10.1109/SWCT.1964.27}{doi:doi.ieeec%
omputersociety.org/10.1109/SWCT.1964.27}, pp.~133--137.

\bibitem{EFH80}
P.~Erd{\H{o}}s, S.~Fajtlowicz, and A.~J. Hoffman, \emph{Maximum degree in
  graphs of diameter {$2$}}, Networks \textbf{10} (1980), no.~1, 87--90.

\bibitem{FP10}
R.~Feria-Pur\'on and G.~Pineda-Villavicencio, \emph{On bipartite graphs of
  defect at most 4},  (2010), submitted.

\bibitem{Hey96}
M.~C. Heydemann, \emph{Cayley graphs and interconnection networks}, {Graph
  Symmetry: Algebraic Methods and Applications} (G.~Hahn and G.~Sabidussi,
  eds.), NATO ASI Series, Series C: Mathematical and Physical Sciences, vol.
  497, Kluwer Academic Publishers, Dordrecht, 1996.

\bibitem{HS60}
A.~J. Hoffman and R.~R. Singleton, \emph{On {Moore} graphs with diameter 2 and
  3}, IBM Journal of Research and Development \textbf{4} (1960), 497--504.

\bibitem{Jo92}
L.~K. J{\o}rgensen, \emph{Diameters of cubic graphs}, Discrete Applied
  Mathematics \textbf{37/38} (1992), 347--351,
  \href{http://dx.doi.org/10.1016/0166-218X(92)90144-Y}{doi:10.1016/0166-218X(%
92)90144-Y}.

\bibitem{KT81}
K.~Kurosawa and S.~Tsujii, \emph{Considerations on diameter of communication
  networks}, Electronics and Communications in Japan \textbf{64A} (1981),
  no.~4, 37--45.

\bibitem{MNP09}
M.~Miller, M.~Nguyen, and G.~Pineda-Villavicencio, \emph{On the nonexistence of
  graphs of diameter 2 and defect 2}, Journal of Combinatorial Mathematics and
  Combinatorial Computing \textbf{71} (2009), 5--20.

\bibitem{PM1}
M.~Miller and G.~Pineda-Villavicencio, \emph{Complete catalogue of graphs of
  maximum degree 3 and defect at most 4}, Discrete Applied Mathematics
  \textbf{157} (2009), no.~13, 2983--2996,
  \href{http://dx.doi.org/10.1016/j.dam.2009.04.021}{doi:10.1016/j.dam.2009.04%
.021}.

\bibitem{MS05}
M.~Miller and R.~Simanjuntak, \emph{Graphs of order two less than the {Moore}
  bound}, Discrete Mathematics \textbf{308} (2008), no.~13, 2810--2821,
  \href{http://dx.doi.org/10.1016/j.disc.2006.06.045}{doi:10.1016/j.disc.2006.%
06.045}.

\bibitem{Mol05}
S.~G. Molodtsov, \emph{General theory of information transfer and
  combinatorics}, Lecture Notes in Computer Science, vol. 4123/2006,
  ch.~Largest Graphs of Diameter 2 and Maximum Degree 6, pp.~853--857, Springer
  Berlin/Heidelberg, 2006.

\bibitem{NMb}
M.~H. Nguyen and M.~Miller, \emph{Structural properties of graphs of diameter 2
  with defect 2}, AKCE International Journal of Graphs and Combinatorics \textbf{7} (2010), no.~1,  29--43.

\bibitem{PM06}{
G.~Pineda-Villavicencio and M.~Miller, \emph{On graphs of maximum degree 5,
  diameter {$D$} and defect 2}, Proceedings of MEMICS 2006, Second Doctoral
  Workshop on Mathematical and Engineering Methods inComputer Science (Mikulov,
  Czech Republic) (L.~Matyska, A.~Ku\v cera, T.~Vojnar, Z.~Kot\'asek, D.~Anto\v
  s, and O.~Kraj\'i\v cek, eds.), Oct 2006, pp.~182--189.
  }
\end{thebibliography}
\end{document}